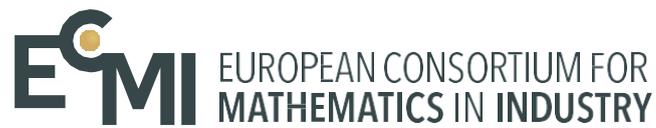

First Virtual ECMI Modeling Week 2020
Peter the Great St.Petersburg Polytechnic University

# BOLTED ASSEMBLY OPTIMIZATION

Margarita Petukhova (Peter the Great St.Petersburg University)

Tahir Miriyev (Università degli studi di Verona)
Yicheng Chen (Beihang University)
Lorenz Huber (Technische Universität Berlin)
Marco Frieden (Universidade de Coimbra)
Zeina Amer (Technische Universität Berlin)

July 2020

# Contents

# 1   Introduction

At the first stage of aircraft assembly, the consistent parts are joined together via installation of temporary fasteners in holes drilled at pre-defined positions. Such a procedure is performed with the aim to minimize the residual gap between parts and prepare them for the following stage of riveting. Although the optimal location of fasteners is computed in advance, it turns out that the order of installation matters as well. As every consecutive fastener is installed, the load in previously installed fasteners can decrease below the acceptable threshold, thereby causing the residual gap to vary uncontrollably. Such a phenomena poses a challenging optimization problem, with a fast and efficient algorithm to be designed to minimize the ultimate residual gap by finding the optimal installation order. Riveting can be performed only once all the fasteners have reached an optimum load force and residual gap is below the acceptable cap.

The work in *Numerical approach for airframe assembly* presents the overview of the ASRP software complex, developed for simulating the installation process and computing quantitative data such as gaps, load forces and so on. This software was utilized by all the algorithms presented in this report.

This report presents several algorithms developed to determine the optimal sequence of actions.

The simulation considers a simplified model of two flat parts with 40 holes drilled at pre-defined XY coordinates and numerated from 0 to 39. The initial uniform gap between parts is 6mm where no fasteners have been installed yet. At each step one of the following actions can be performed: a new fastener can be installed at an empty hole or an existing fastener can be tightened, the latter of which we will define as refastening. In both cases the force applied at fastener installation is set to 1000 Newton, initially equal to the load force value of an installed fastener, subject to change as each consecutive fastener is installed.

As a consequence of installing or tightening a fastener, the parts are displaced in relation to one another. The displacement depends on load forces, the structural rigidity of the parts and relative distances between holes. Optimal load force is achieved, when all fasteners apply a force of $1000 \pm 25$ Newton. Once all installed fasteners have an optimal forces, there is no need for refastening and pre-riveting stage is completed. At that point the fasteners at the selected set of holes apply an optimal force to the parts and for the chosen set of holes a minimal gap is achieved.

The ASRP software complex takes as an input a sequence of numbers 0 to 39, corresponding to holes indexing. In sequential order, a fastener is installed at the hole or an existing fastener is refastened, depending on if the number appears for the first time or is being repeated. After every action the software returns the resulting gaps and forces in all holes.

The main objective of developed algorithms is to achieve optimal force load with the minimal number of actions performed. With an ideal number and distribution of fasteners the gap in the optimal state is defined by mean = 0.038 and standard deviation = 0.042. These values are achieved by simulating a simultaneous installation of 20 fasteners. The optimal configuration of 20 fasteners has been provided to us. However, note that under real circumstances it is not possible to install all fasteners at the same time, hence these statistics are abstract. Hence the goal is to get as close as possible to an uniform distribution.

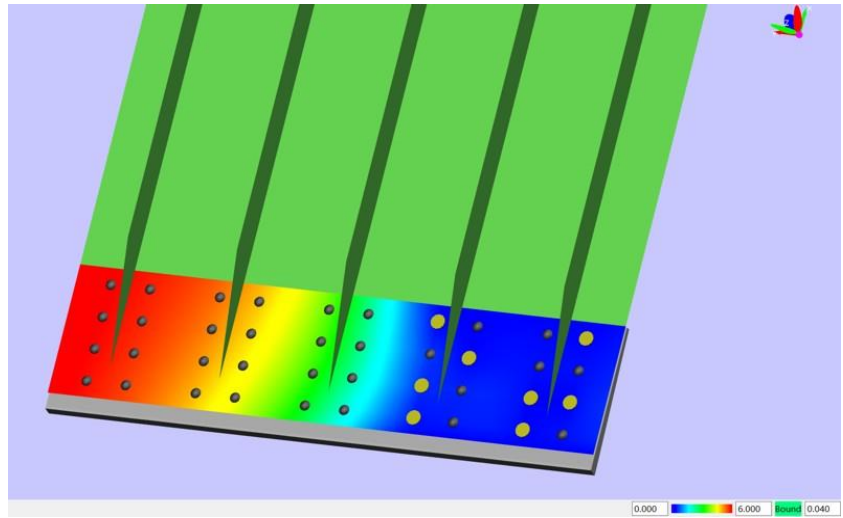

Figure 1: Simulated parts as viewed from above with 8 fasteners installed and a colored gradation indicating the gap, red (large gap) to blue (small gap).

## 2 Report content

### 2.1 Mathematical Preliminaries

We consider a Finite Element Model (FEM) representing the assembly process. The unknowns are displacements and rotations in all FEM nodes that form the vector $U = (U_i)_{i=1}^{N_{DOF}}$ where $N_{DOF}$ is the total number of Degrees Of Freedom (DOFs) included in the model. The optimization problem can be formulated in terms of Quadratic Programming as finding a minimum of energy functional:

$$\min \frac{1}{2} U^t \cdot K \cdot U - F^t \cdot U$$

subject to

$$(A \cdot U)_i < g_i \quad i = 1, ..., N$$

where the admissible class of displacements is $A_C = \{U | U = U_0 \text{ in } \Gamma, A \cdot U \leq G\}$,

A = linear operator that defines constrained DOFs in contact panels,

$F = (F_i)_{i=1}^{NDOF}$ the vector of force loads,

K = fastener stiffnes matrix

$G = (g_i)_{i=1}^{NDOF}$ the vector of initial gaps defining the constraints of displacements between corresponding nodes

The original problem is considered for FEM with thousands of nodes, which is computationally expensive. Instead, we solve the problem belonging to denoted class R we can reduce its dimension by eliminating all unknowns outside the junction area. The junction area is the area of "computational" nodes, where about 20% of total nodes are located. The process of joining panels occurs on the junction area.

For that we choose the system of 40 computational nodes $CN = (cn_1, cn_2, cn_3, ..., cn_{40})$ in junction area, for example, as in above figures. Then, we build the reduced stiffness matrix $K_c$. Reduced rigidity matrix characterizes the response of computational nodes to the applied load. To be more precise, $K_C^{-1} = \{k_{ij}\}_{i,j=1}^{40}$ is inverse matrix to $K_c$ and $k_{ij}$ is a displacement of $j^{th}$ node caused by the unit load applied to $i^{th}$ computational node. Let us now consider a method for computing matrix.

We divide whole displacement vector into two parts: $U = \begin{pmatrix} U_c \\ U_r \end{pmatrix}$ where $U_c$ is a vector of displacements in junction area, $U_r$ is vector of displacements in all other nodes. Then global stiffness matrix K can be rewritten in blocks as:

$$\begin{pmatrix} K_{cc} & K_{cr} \\ K_{cr}^t & K_{rr} \end{pmatrix} \cdot \begin{pmatrix} U_c \\ U_r \end{pmatrix} = \begin{pmatrix} F_c \\ F_r \end{pmatrix}$$

Calculating Schur complement for the block $K_{rr}$ we obtain reduced rigidity matrix $K_c$ instead of block $K_{cc}$. Also this matrix can be computed using the formula:

$$K_c = K_{cc} - K_{cr} \cdot K_{rr}^{-1} \cdot K_{cr}^t$$

Now we can formulate equivalent minimization problem for UC:

$$\min \frac{1}{2} U_c^t \cdot K_c \cdot U_c - F_c^t \cdot U_c$$

where $A_c = \{U_c | A \cdot U_c \leq G\}$

## 2.2 MaxGap algorithm

As a rather simple and very intuitive approach we implemented this algorithm to try to solve our problem. The general idea behind it is to always tighten the fastener that has the biggest gap. The following points may help to understand the approach.

- The first fastener is chosen at random or by the user. Then the algorithm begins:
    (1) Check which hole has the biggest gap value.
    (2) Tighten the fastener in the hole found in (1).
- Repeat steps (1) and (2) until the chosen stopping criteria are met.

## 2.3 MaxMinDivide algorithm

### 2.3.1 Description

Similar to the algorithm presented before, we wanted to tryout different approaches and get better results and implemented the MaxMinDivide algorithm. The intuitive idea behind this algorithm is to set the fasteners gradually from an already tight hole to the loosest, the one with the biggest gap. The steps it follows are presented in the following points:

(1) Install the first fastener. This choice is given by the user.
(2) Search for the biggest gap in all chosen holes, the ones with and without fasteners.
(3) Search for the smallest gap in all chosen holes, the ones with and without fasteners.
(4) Calculate the difference between the maximum and the minimum and divide it by n. This n is chosen as wished, as long as it is positive.
(5) Search which hole is closest to the sum of the values of 3. and 4..
(6) The hole returned in 5. is where the program will either mount a fastener or tighten an already installed one.
(7) Repeat from point 2. on until the chosen stopping criteria have been satisfied or the maximum amount of iterations has been reached.

### 2.3.2 Simulation

In this work we analyse for $n$ implemented as 2, 4 and number of fasteners to be used (#fasteners). We chose to perform the same tests on each:

- Starting at hole 10, put fasteners into holes 10, 23 and 25.
- Starting at hole 24, put fasteners into holes 3, 9, 22, 24 and 36.
- Starting at hole 19, put fasteners into the 20 holes that would give the best result, could we tighten 20 the fasteners all at the same time.

In this report we will analyse the latter choice, with 20 fasteners. The limit for the number of simulations did vary depending on how the results were developing and the time needed to compute them.

### 2.3.3 Code

https://drive.google.com/drive/folders/1g2ppr2Q0byNOWhI_8J0qbfsiERKf6-YH?usp=sharing

## 2.4 Blockwise search

### 2.4.1 Description

The main idea in this algorithm is to make use of the geometry of the panel. As you can see in figure 1, the 40 holes are evenly distributed in 5 blocks of 8 nodes each. By limiting the selection of the next screw placement to a specific environment, namely the neighboring block, we will be decreasing the gaps in a uniform manner regarding the panel area.
The choice of placement of the next fastener is made as follows:

1. Calculate which block is the nearest neighboring block B to the one with the last installed fastener.
2. From said neighboring block B, determine the smallest gap $s$ and the largest gap $l$ and set avg = $\frac{s+l}{2}$
3. Determine the fastener/hole with the closest gap to avg-value in the block B.
4. Do the installment in that determined hole, being whether refastening the fastener already there or installing a new one in.

### 2.4.2 Simulation

Due to processor inconvenience from my side, the simulation could unfortunately not be run until it was done. We had to terminate manually after the 30th step, which was sufficient to see that the gap and force means converge to the desired values, which you can see in the figure 3.3 in the discussion.

### 2.4.3 Code

https://drive.google.com/drive/folders/11XnGjRxAhYXV8dD1LkC6IS1hZLTdluvP?usp=sharing

## 2.5 MAXPERIM algorithm

### 2.5.1 Description

Before building the logic behind MAXPERIM, a few simple scenarios were simulated and the following crucial observations and decisions were made:

- Firstly, it was important to form a criteria on choosing the next node in a way that it minimally affects the gap in the already installed fasteners. From simulations we noticed direct correlation between $d(n_j, n_{j+1}) \cong gap(n_j)$, i.e. distance to the next node and the change of gap value in the previous node(s). It became evident that the larger distances are, the smaller gap changes are. This observations lies in the core of the MAXPERIM algorithm.

- Tolerance ranges were defined as 0.01 and 0.02 for final mean and final stdev respectively. After every action the algorithm re-calculates current statistics, and once they lie within tolerance intervals ($currentmean \pm meantol$), ($currentstdev \pm stdev\ tol$), the algorithm stops.

- Another important criteria is the action decision criteria: whenever $Force(n_j) < 990N$ we tighten a fastener. Otherwise we install a new fastener.

Given all above-mentioned criteria, the logic behind MAXPERIM is straightforward. For N fasteners already installed, the next fasteners is chosen in a way that maximizes the area of a polygon formed by these fasteners. After running the algorithm, we reached the optimal gap distribution in 36 actions. The results are displayed in Figure 2. 2

### 2.5.2 Simulation

| ACTIONS | FASTENING ORDER |
|---|---|
| new fastener' | 1 |
| new fastener | 38 |
| new fastener | 3 |
| tightened fastener | 1 |
| new fastener | 2 |
| new fastener | 8 |
| new fastener | 6 |
| new fastener | 11 |
| new fastener | 10 |
| tightened fastener | 3 |
| tightened fastener | 38 |
| tightened fastener | 3 |
| new fastener | 16 |
| new fastener | 14 |
| new fastener | 19 |
| tightened fastener | 11 |
| tightened fastener | 1 |
| tightened fastener | 11 |
| tightened fastener | 2 |
| tightened fastener | 11 |
| tightened fastener | 8 |
| tightened fastener | 11 |
| tightened fastener | 6 |
| tightened fastener | 11 |
| new fastener | 18 |
| new fastener | 24 |
| new fastener | 22 |
| new fastener | 27 |
| new fastener | 26 |
| new fastener | 32 |
| tightened fastener | 38 |
| tightened fastener | 10 |
| tightened fastener | 38 |
| new fastener | 30 |
| new fastener | 34 |
| new fastener | 35 |

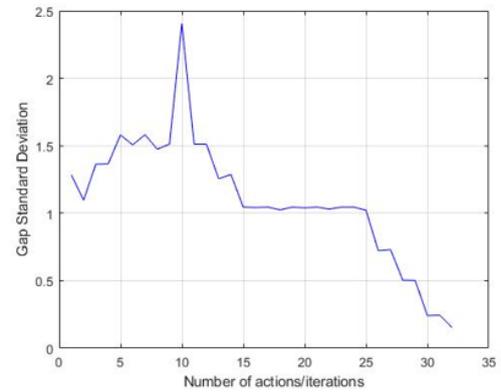
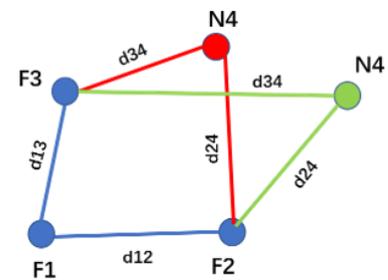
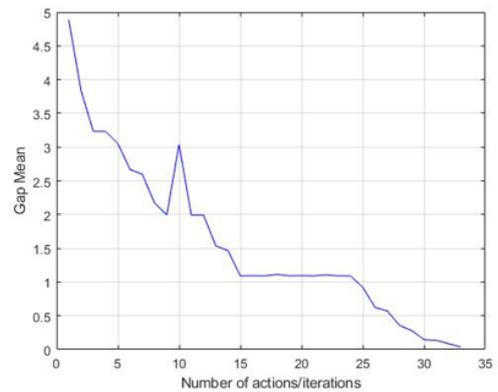

Figure 2: Performance of the MAXPERIM algorithm

### 2.5.3 Code

https://drive.google.com/drive/folders/1cIU40GjbzaAffebv4vXO2ofIz0gyAh0X?usp=sharing

## 2.6 Gap gradient-based search

### 2.6.1 Description

The goal is to minimize the global gap while ensuring that the gaps at each hole do not vary too much, namely, we want to both minimize the mean of gaps and the variance of gaps. Consequently, consider the objective function to minimize as:

$$L = \lambda \sigma^2 + (1 - \lambda) \cdot \frac{1}{N} \sum_{i=0}^{N-1} g_i \qquad (2.1)$$

where $g_i$ stands for the gap at the *ith* hole. $\sigma^2$ is the variance of all gaps. $N$ is the total number of holes. The constant $\lambda$ is used for balancing the weights between the mean of gaps and the variance of gaps.

A special point of this optimization problem is that the sequence does not explicitly express object function L. Then the optimal sequence should make the biggest contribution to this objective: $min(L)$

### 2.6.2 Updating method

This method uses the difference of gaps at each step to determine the location of the next fastener.

Assume that there are $k$ fasteners already installed, the location of the $k+1$ fastener can be chosen by the following approach:

(1) Calculate the gaps when there are $k$ fasteners installed, get the gaps vector $G_k$ at each hole.

(2) Subtract $G_k$ and the former gaps vector $G_{k-1}$, get the difference vector $D_k$, each element in $D_k$ is $d_i^k (i = 0, 2, 3, \cdots, N - 1)$.

(3) Assume the indices of all holes are set A, the indices of the $M$ fasteners are set B, and the indices of the $k$ fasteners are set E. Choose the index $i$ from the set $A \cap (B - E)$ where the corresponding $_i d^k$ is the largest in $D_k$. Then install the fastener at the hole *ith* and add index $i$ to set E.

(4) Repeat (1)-(3) until all the $M$ fasteners are installed.

The location of the first fastener can be obtained by running the algorithm starting with each of the $M$ fasteners and find out which one minimizes the objective function.

### 2.6.3 Simulation

The simulation parameters are set as $N = 40$, $\lambda = 0.6$, $M = 20$. The indices of these 20 fasteners are [1 2 3 6 8 10 11 14 16 18 19 22 24 26 27 30 32 34 35 38].

Fig. 3 shows the changes in gaps during the installation of 20 fasteners. Based on the algorithm proposed, the gaps decrease sharply. We have also noticed these facts:

(1) The gaps at each hole seldom rise during the whole procedure.

(2) Generally, the proposed algorithm tends to choose the next hole close to the current hole.

(3) After 20 steps, all the gaps are relatively much smaller than they were initially. And there is no considerable variance between them.

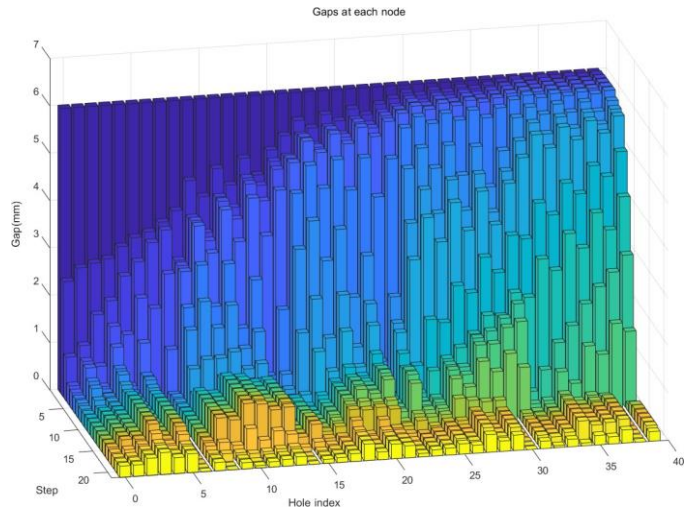

Figure 3: The changes in gaps during the installation of 20 fasteners (starting from hole 3)

Fig. 4 shows the change of loss (the value of the objective function) versus steps. We can notice that the loss rises in the first few steps and decreases sharply in the middle, then changes little until the end. The rise of the loss at the beginning is reasonable: Initially, the variance of the gaps is zero, when there are only a few fasteners installed, the mean of gaps does not change dramatically, while the variance of gaps is greatly lifted, which causes the loss to rise in the first few steps.

Fig. 5 shows the mean of gaps, variance of gaps, and loss after 20 steps versus the index of starter hole. Based on fig. 5, we can confirm that the optimal starter hole is the $10th$ hole.

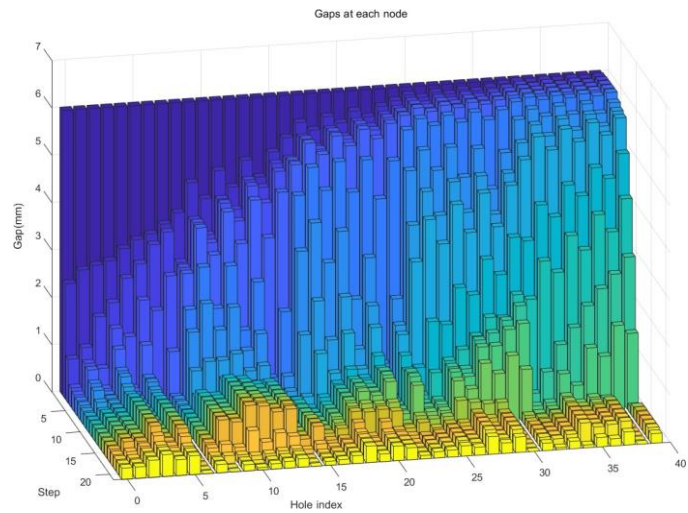

Figure 4: Loss versus steps (start from hole 3)

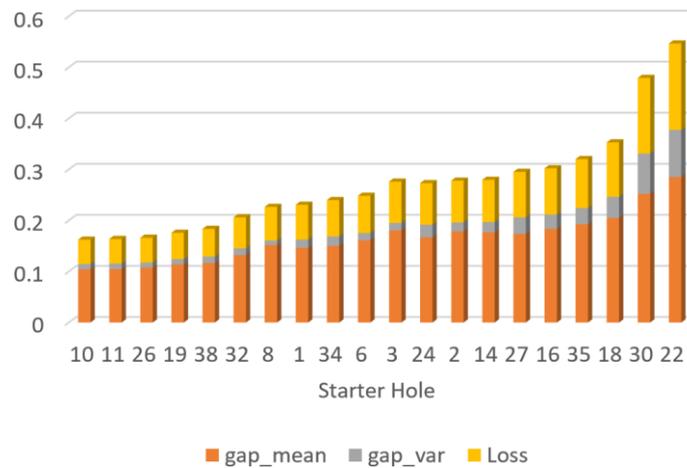

Figure 5: The mean of gaps, the variance of gaps, and loss after 20 steps versus the index of starter hole.

### 2.6.4 Code

https://github.com/Amos-Chen98/Gap-gradient-based-search_ECMI-Modelling
-Week2020_Q1

## 2.7 KF algorithm

### 2.7.1 Description

The research in *Numerical approach for airframe assembly* is centered on a quadratic optimization problem with constraints. The objective function is

$$E(U) = \frac{1}{2} U^T \cdot K \cdot U - F^T \cdot U, \tag{2.2}$$

where $E(U)$ describes the contact stress between the parts as the force $F$ is applied to the fasteners. The global stiffness matrix $K$ represents the structural rigidity of the parts. The energy functional $E$ is minimized over the displacement $U$, to reach the optimal gap between two parts. The constraints are not discussed here.

We would like to determine the next fastener, to be installed or refastened, with the greatest influence on displacement. To determine which action will result in the largest displacement, iteratively every fastener is fastened, the resulting displacement is recorded, and the fastener is released to its original force. All possible displacements are computed, compared and the action resulting in the largest displacement is chosen.

### 2.7.2 Simulation

The algorithm solves the quadratic objective function without constraints, i.e. the constraints are ignored. The global stiffness matrix is positive definite, thus a unique solution can be determined by solving a system of linear equations

$$K \cdot U = F. \tag{2.3}$$

An initial displacement vector $U_{init}$ is calculated with the forces currently applied to the fasteners $F_{init}$.

The fasteners are indexed $i = 0, ..., n$. In every iteration the force of a fastener $j$ is set to the optimal tension $f_j = 1000$, $F_j = [f_i]$, and the resulting $U_j$ is calculated. The displacements are compared by determining the Euclidean norm of the absolute difference $U_{init} - U_j$.

### 2.7.3 Code

https://drive.google.com/drive/folders/1IdiuDrd0Ws4NwNcs-ZZQb9WUBZ4_tjHC?usp=sharing

# 3 Discussion and conclusions

## 3.1 MaxGap

This first algorithm, which we comically named 'the bicycle', due to it's simplicity, did surprisingly well. Having in mind that we wanted to maximise the applied force and minimise the gap at each fastener it was able to tighten all 20 fasteners in 46 iterations at $1000 \pm 25$ Newton.

## 3.2 MaxMinGap

This second implementation was able to get some good results as well but did not converge as the more simpler MaxGap algorithm as even after 60 iterations, neither of the 3 algorithms was able to converge to all fasteners having $1000 \pm 25$ Newton in terms of force:

| $n = \#fasteners$ | Gap Var | Gap Mean | Force Var | Force Mean |
|---|---|---|---|---|
| 20 it | 1.0350 | 0.8444 | 145150 | 204.050 |
| 30 it | 0.0099 | 0.1203 | 197590 | 353.825 |
| 60 it | 0.0021 | 0.0458 | 244810 | 462.275 |

| $n = 2$ | Gap Var | Gap Mean | Force Var | Force Mean |
|---|---|---|---|---|
| 20 it | 1.0182 | 0.1995 | 168050 | 277.025 |
| 30 it | 0.0092 | 0.1459 | 200500 | 303.500 |
| 60 it | 0.0227 | 0.0837 | 256560 | 475.625 |

| $n = 4$ | Gap Var | Gap Mean | Force Var | Force Mean |
|---|---|---|---|---|
| 20 it | 0.1041 | 0.4559 | 158590 | 225.325 |
| 30 it | 0.0040 | 0.0957 | 194310 | 349.825 |
| 60 it | 0.0036 | 0.0543 | 256730 | 475.775 |

## 3.3 Blockwise search

The algorithm tends to prioritise refastening already installed fasteners, which probably is not the ideal approach. Yet you can see in the figure below, that the gap and force means are converging to the beforehand calculated desired values.

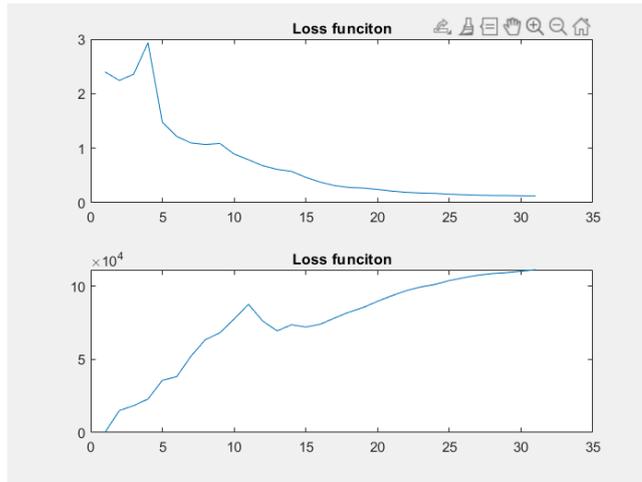

Figure 6: Top: gap mean, Bottom: force mean

## 3.4 MAXPERIM algorithm

Indeed, MAXPERIM was not the very first algorithm proposed by the author. Due to the deficiency of MATLAB's built-in functions capable of calculating the area of any polygon (given the vertices in no particular order), the author could not realize the initial algorithm he thought of, MAXAREA. The principle behind two algorithms is slightly different, however the MAXAREA is expected to produce better results. As can be guessed from the name, MAXAREA tries to maximize the area (instead of a perimeter) of a polygon formed by previously installed fasteners and the next fastener. At times the choice of the next node decided by MAXAREA and MAXPERIM will overlap, but not always. MAXAREA could avoid non-steady decline of variance and mean.

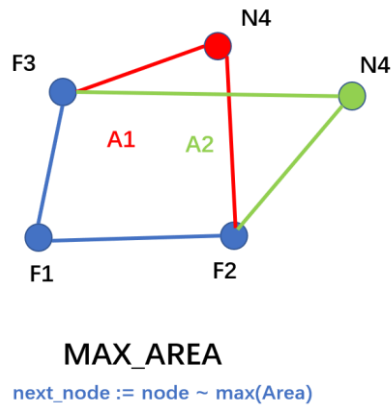

Figure 7: MAXAREA algorithm illustrated for a simple case of 3 installed fasteners and the 4th one to be chosen red or green)

In 36 actions/steps the algorithm succeeded to reach the ideal mean and standard deviation. Since there was no criteria on the minimal number of actions an algorithm should take, it is hard to tell whether the outcome of MAXPERIM is optimal or not. It should be noted that the posed research question defines a totally new field of research, therefore making it hard to define any sensible criteria. The author was free to state his criteria, in terms of *mean tolerance* = 0.01, *stdev tolerance* = 0.02 and *force tolerance* = 10$N$, these values however should be negotiated with experienced engineers and scientists in AIRBUS. From Figure and one can't observe a steady decline of both mean and st.dev values. Avoiding spikes at particular time instants could be another criteria, for example. MAXAREA seems to be capable of producing better results in that sense.

## 3.5 Gap gradient-based search

Based on all the simulations under the proposed algorithm, the gap gradient-based search method has been tested thoroughly, it performs well in the assembly scenario of choosing the best sequence of installation with fixed steps and fixed location. There are remarkable advantages of this algorithm:

(1) Based on backward differential, it is neat in logic and easy to follow in the real industrial production.

(2) Time complexity is linearly related to the number of nodes, which makes the algorithm remain efficient even the number of nodes is extremely huge.

(3) The performance of the algorithm is relatively high since it can reduce the mean of gaps to $0.104 mm$ and the variance of gaps to $0.010 mm^2$ within 20 steps.

There are still some disadvantages of this algorithm, for example, the selection of the first node depends on experiments, which might not be practical in some certain situations.

We can come to this summary based on the proposed algorithm:

The best node to start with is the $10th$ node. The best sequence is [10 8 1 2 6 14 22 30 38 35 32 24 26 27 34 3 11 16 18 19], subject to $N = 40$, $\lambda = 0.6$, $M = 20$.

## 3.6 KF algorithm

The algorithm attempts to predict the displacement of the parts without having to execute expensive calculations. Instead of having to minimize a quadratic problem with constraints, for every prediction a system of linear equation is solved. Thus it is possible to make more predictions.

Figure 8 shows the forces applied to the fasteners and the resulting gaps present at the fasteners after 7 actions and 55 actions. The algorithm starts off by installing and repeatedly fastening a centrally located set of fasteners. It then progresses to the periphery, at which point fasteners can be installed without having to be refastened repeatedly.

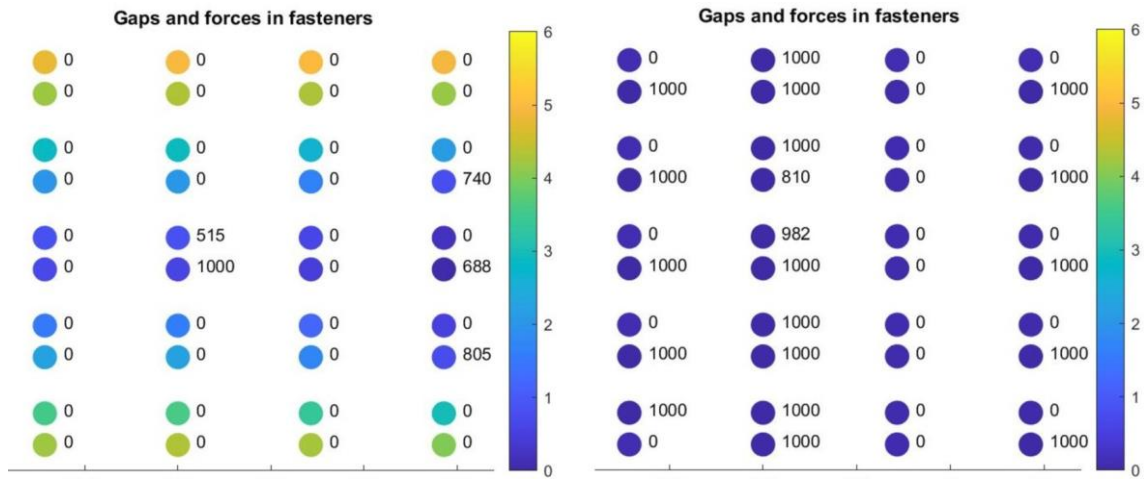

Figure 8: Gaps and forces in fasteners after 7 actions (left) and 55 actions (right). Every dot stands for a fastener, the locations of the dots in the graph represent the spacial distribution of the fasteners on the parts to be fastened. The color shows the gap between the parts in *mm*. The adjacent number indicates the force applied.

The algorithm was able to fasten all fasteners of the 20 fastener test scenario to a force of 1000 +/-25 Newton in 56 steps: 19 18 30 14 22 19 18 11 2 16 27 8 10 3 35 1 32 34 38 24 26 6 27 30 27 14 27 22 27 19 27 18 27 11 27 2 27 16 27 35 8 35 10 35 3 35 34 1 34 32 34 26 38 26 24 26.

Disregarding the constraints of the optimization problem will not accurately calculate the displacement of the parts. The approach is based on the assumption, that the simplified calculation might serve as an indicator for which fastener to fasten next.

The algorithm is cheap to execute and is able to fasten the fasteners as required, still doubt remains if the predictions have a strong connection to reality. This approach seems to confirm the idiom "*One can not have a cake and eat it.*"

# 4  Group work dynamics

The group work was organized productively. Despite the fact that none of us specialized in algorithm design and had limited background in the field of airplane assembly, each of us individually managed to propose interesting ideas on how to design potentially well-performing algorithms.

We began research by studying the papers provided by our instructors. Once we have got more complete understanding of underlined mathematics and physics, we formulated a few assumptions that could lie in the core of our algorithms. After testing those assumptions through software simulations, each team member came up with a creative idea and decided to focus on his/her algorithm, while regularly updating the others on progress.

Everyone was cooperative and supportive. Whether it was about writing a code or simulating assumptions, preparing the final presentation or understanding theoretic material, each team member willingly dedicated his/her time to help. We believe that helped each of us to succeed in our endeavors, both individually and as a team.

We are utterly thankful to our instructors, Dr. Maria Churilova and Dr. Margarita Petukhova, whom we would always approach if we had some theoretic questions or technical problems. Overall, the atmosphere was cooperative, motivating and very friendly.

# 5 Instructor's assessment

The problem posed for the group is brand new and challenging even for me who works in this field for more than 10 years. Students had to study the mathematical basics of the problem, to manage with the code that was provided for them and to develop their own algorithms. It was a pleasure to see how thoroughly they approached the issue, each of them could propose and realize the own solution.

The group was very well organized from the first to the last day of the Modelling Week. I would like to thank them for the real teamwork: being apart, they supported each other as if they were together. It was the first experience of online conference for all of us, and I really appreciate what they have done.

# 6 REFERENCES

- Petukhova M., Lupuleac S., Shinder Y., Smirnov A., Yakunin S.,Bretagnol B. (2014).Numerical approach for airframe assembly. Journal of Mathematics for Industry.